\documentclass{amsart}
\usepackage{amsmath}
\usepackage{amssymb}
\usepackage{amsbsy}
\usepackage{mathrsfs}
\usepackage[all]{xypic}

\newtheorem{theorem}{Theorem}[section]
\newtheorem{dummy}[theorem]{}
\newtheorem{prop}[theorem]{Proposition}
\newtheorem{lemma}[theorem]{Lemma}
\newtheorem{cor}[theorem]{Corollary}

\newtheorem{example}[theorem]{Example}
\newtheorem{remark}[theorem]{Remark}

\numberwithin{equation}{section}

{
\vspace{0.25cm}
}
{
\vspace{1cm}
}

{
\vspace{0.25cm}
}
{
\vspace{1cm}
}

\begin{document}

\title{Ossa's Theorem via the Kunneth formula}

\author{Robert Bruner, Khairia Mira, Laura Stanley and Victor Snaith}

\date{28 June 2010 
} 

\begin{abstract} 

Let $p$ be a prime. We calculate the  connective unitary K-theory of the smash product of two copies of the classifying space for the cyclic group of order $p$, using a K\"{u}nneth formula short exact sequence. As a corollary, using the Bott exact sequence and the mod $2$ Hurewicz homomorphism we calculate the connective orthogonal K-theory of the smash product of two copies of the classifying space for the cyclic group of order two.
\end{abstract}

\maketitle

\section{Introduction}

This paper arose as a result of discussions during a graduate course at the University of Sheffield during 2008. In order to introduce Frank Adams' technique of constructing homology resolutions as realisations of iterated cofibrations of spectra a simpler example than the classical Adams spectral sequence was needed. We had the spectrum $bu$ to hand but, in order to postpone the algebraic intricacies of spectral sequences, what was required was an example whose geometric resolution gave rise to a short exact sequence rather than a spectral sequence. As it happens the $bu$-resolution of $B{\mathbb Z}/2$ yields such an example, which was simple enough for the purposes of the course. At that point, John Greenlees mentioned the existence of \cite{Os89}, which prompted the writing of \S2. In \S2 we use the $bu$-resolution of \S\ref{2.1} to calculate $bu_{*}(    B{\mathbb Z}/p \wedge   B{\mathbb Z}/p)$ in terms of $bu_{*}(    B{\mathbb Z}/p )$ and mod $p$ Eilenberg-MacLane spaces (Theorem \ref{2.12}).

We shall merely compute connective K-groups. The papers \cite{Br99},  \cite{JW97}  and \cite{Os89} derive equivalences of spectra involving  $bu \wedge  B{\mathbb Z}/p \wedge   B{\mathbb Z}/p$ and $bo \wedge  B{\mathbb Z}/p \wedge   B{\mathbb Z}/p$ which imply our results upon taking homotopy groups. Importantly, unlike  \cite{Br99},  \cite{JW97}  and \cite{Os89}, we do not resort to Adams spectral sequences to construct the essential algebraic homomorphism
\[  \tilde{\mu}_{*} : bu_{*}(    B{\mathbb Z}/p \wedge   B{\mathbb Z}/p)  \longrightarrow   \vee_{i=1}^{p-1}  \ bu_{*}(  S^{2i}  \wedge   B{\mathbb Z}/p)   \]
of \S\ref{2.11}. 

Our method does not yield a homomorphism $\tilde{\mu}_{*} $ induced on homotopy from a map of spectra, it is merely an algebraic homomorphism, as explained in \S\ref{2.13}. However, since $\tilde{\mu}_{*} $ is closely related to the map induced by the multiplication on $B{\mathbb Z}/p$ it is virtually invariant under switching the  $B{\mathbb Z}/p$-factors, which may prove useful in calculations of $bu_{*}$ of other $p$-groups.

In \cite{Br99} and \cite{Os89} $bo_{*}$-analogues of the $bu$-result are offered when $p=2$. Our $bo_{*}$ calculations are consistent with the results proved in \cite{Br99} and highlight the errors in the $bo$-analogue asserted in \cite{Os89}.
Consider the cofibration
\[  \Sigma bo \stackrel{\eta}{\longrightarrow}  bo   \stackrel{c}{\longrightarrow} bu  ,  \]
discovered by Raoul Bott during the proof of his famous Periodicity Theorem. Smashing this with $X$ and taking homotopy groups yields the Bott sequence for $X$. In \S3
we compute $bo_{*}( B{\mathbb Z}/2 \wedge  B{\mathbb Z}/2)$ by comparing the Bott sequence for $B{\mathbb Z}/2 \wedge  B{\mathbb Z}/2$ with that for $B{\mathbb Z}/2 $ and with mod $2$ homology. Our calculations are relevant to 
\cite{BG03} and \cite{Arf10}, for example.

\section{The connective unitary case}
\begin{dummy}
\label{2.1}
\begin{em}

Let $bu_{*}$ denote connective unitary K-homology on the stable homotopy category of CW spectra \cite{Ad74} so that if $X$ is a space without a basepoint its unreduced $bu$-homology is $bu_{*}( \Sigma^{\infty}X_{+})$, the homology of the suspension spectrum of the disoint union of $X$ with a base-point. In particular $bu_{*}( \Sigma^{\infty}S^{0})= {\mathbb Z}[u]$ where ${\rm deg}(u) = 2$. Let $p$ be a prime and consider the cofibration of pointed spaces
\[      B{\mathbb Z}/p  \stackrel{i}{\longrightarrow} BS^{1} \stackrel{\pi}{\longrightarrow}  W_{p}   \]
where $i$ is induced by the inclusion of the cyclic group of order $p$ into the circle. This cofibration maps to the fibration
\[     B{\mathbb Z}/p  \stackrel{i}{\longrightarrow} BS^{1} \stackrel{Bp}{\longrightarrow}  BS^{1}  \]
and the comparison for mod $p$ and integral unreduced singular homology yields the following result:
\end{em}
\end{dummy}
\begin{lemma}{$_{}$}
\label{2.2}
\begin{em}

For all $j$, $H_{j}(W_{p} ; {\mathbb Z}) \cong H_{j}(BS^{1}; {\mathbb Z})$  being ${\mathbb Z}$ when $j \geq 0$ is even and zero otherwise.
\end{em}
\end{lemma}

From the Atiyah-Hirzebruch spectral sequence (\cite{Ad74} p.47) we obtain the following result, which also follows from the Thom isomorphism $bu_* (W_p) \cong bu_*( BS^1)$, since $W_p$ is Thom complex of the $p$-th tensor power of the canonical complex line bundle, by \S\ref{2.1}. 
\begin{cor}{$_{}$}
\label{2.3}
\begin{em}

Both  $bu_{*}(\Sigma^{\infty} W_{p})$ and $bu_{*}(\Sigma^{\infty} BS^{1})$ are free modules over $bu_{*}( \Sigma^{\infty}S^{0}) = {\mathbb Z}[u]$.
\end{em}
\end{cor}
\begin{dummy}
\label{2.4}
\begin{em}
The Atiyah-Hirzebruch spectral sequences for for $bu_{*}$ and $KU_{*}$ of $ \Sigma^{\infty} B{\mathbb Z}/p $
both collapse for dimensional reasons and the map between them is injective so that $bu_{*}(\Sigma^{\infty} B{\mathbb Z}/p)$ injects into $KU_{*}(\Sigma^{\infty} B{\mathbb Z}/p)$ which, by the universal coefficient theorem for $KU$-theory \cite{At62} and the calculations of \cite{At68}, is given by $KU_{2j+1}(\Sigma^{\infty} B{\mathbb Z}/p) \cong  \oplus_{j=1}^{p-1} \ {\mathbb Z}/p^{\infty}$ (\cite{Os89} \S2; see also  \cite{Ka81} Chapter I, \S2) and is zero in even dimensions. 

When $p$ is odd it will be convenient to replace $bu$ by $bu{\mathbb Z}_{p}$, connective unitary K-theory with $p$-adic integers coefficients and similarly for $KU{\mathbb Z}_{p}$. These $p$-adic spectra possess Adams decompositions \cite{Ad69} (see also \cite{Ka81})
\[   bu{\mathbb Z}_{p} \simeq   \vee_{i=1}^{p-1}  \  \Sigma^{2i-2} lu    \  {\rm and}  \   KU{\mathbb Z}_{p}   \simeq   \vee_{i=1}^{p-1}  \  \Sigma^{2i-2} LU   \]
where $lu_{*}(\Sigma^{\infty} S^{0}) \cong {\mathbb Z}_{p}[v]$ where ${\rm deg}(v) = 2p-2$ corresponds to $u^{p-1}$ and multiplication by $u$ translates the summand $\Sigma^{2i-2} lu $ to $\Sigma^{2i} lu $ for $0 \leq i \leq p-2$ and $\Sigma^{2p-4} lu $ to $lu $. $LU$-theory is obtained from $lu$ by localising to invert $v$. In addition there are canonical isomorphisms
\[   bu_{*}( \Sigma^{\infty} B{\mathbb Z}/p) \cong  (bu{\mathbb Z}_{p})_{*}( \Sigma^{\infty} B{\mathbb Z}/p)  \  {\rm and}  \ 
 KU_{*}( \Sigma^{\infty} B{\mathbb Z}/p) \cong  (KU{\mathbb Z}_{p})_{*}( \Sigma^{\infty} B{\mathbb Z}/p)   \]
with $LU_{2j+1}(\Sigma^{\infty} B{\mathbb Z}/p) \cong  {\mathbb Z}/p^{\infty}$.

\end{em}
\end{dummy}
\begin{cor}{$_{}$}
\label{2.5}
\begin{em}
Let $p$ be a prime and let $lu_{*}$ be as in \S\ref{2.4} when $p$ is odd or $lu = (bu{\mathbb Z}_{2})_{*}$ when $p=2$. Then, as a ${\mathbb Z}_{p}[v]$-module, where ${\rm deg}(v_{2i-1}) = 2i-1$,
\[    lu_{*}( \Sigma^{\infty} B{\mathbb Z}/p) \cong  \frac{{\mathbb Z}_{p}[u] \langle  v_{1}, v_{3}, v_{5}, \ldots  \rangle}{(pv_{1} ,  pv_{3}, \ldots , pv_{2p-3}, vv_{2i-1} - pv_{2(p-1) + 2i-1})} . \]
\end{em}
\end{cor}
\vspace{2pt}

{\bf Proof}
\vspace{2pt}

The injection mentioned in \S\ref{2.4} maps $lu_{2i-1}( \Sigma^{\infty} B{\mathbb Z}/p)$ into 
\linebreak
$LU_{2i-1}( \Sigma^{\infty} B{\mathbb Z}/p) \cong {\mathbb Z}/p^{\infty}$. Therefore this group must be cyclic and an order-count in the collapsed Atiyah-Hirzebruch spectral sequence shows that the non-zero groups $lu_{2k(p-1) +2i-1}(\Sigma^{\infty} B{\mathbb Z}/p) \cong {\mathbb Z}/p^{k+1}$ for $i = 1, \ldots , p-1$, generated by $v_{2k(p-1) + 2i-1}$. In $KU_{2i+1}(\Sigma^{\infty} B{\mathbb Z}/p) $ the element $vv_{2k(p-1) +2i-1}$ has order $p^{k+1}$, by Bott periodicity, so we may choose $v_{2(k+1)(p-1) +2i-1}$ so that  $vv_{2k(p-1) +2i-1}= p v_{2(k+1)(p-1) +2i-1}$.      $\Box$
\begin{cor}{$_{}$}
\label{2.6}
\begin{em}

The cofibration of \S\ref{2.1} gives a free ${\mathbb Z}[u] $-module resolution 
\[   0 \longrightarrow  bu_{*}(\Sigma^{\infty} BS^{1})  \stackrel{\pi_{*}}{\longrightarrow}  bu_{*}(\Sigma^{\infty} W_{p})  \longrightarrow   bu_{*-1}(\Sigma^{\infty} B{\mathbb Z}/p)   \longrightarrow 0    \]
as well as similar resolutions for $bu{\mathbb Z}_{p}$ and $lu$.
\end{em}
\end{cor}
\begin{dummy}
\label{2.7}
\begin{em}

If $A$ is a ${\mathbb Z}$-graded group we write $A[n]$ for the graded group with $A[n]_{j} = A_{j+n}$ so that
$bu_{*-1}(\Sigma^{\infty} B{\mathbb Z}/p)$ equals $bu_{*}(\Sigma^{\infty} B{\mathbb Z}/p)[-1]$. By a cell-by-cell induction, for all CW spectra of finite type $X$ the external product gives isomorphisms
\[   \begin{array}{c} 
bu_{*}(\Sigma^{\infty} BS^{1}) \otimes_{{\mathbb Z}[u] }  bu_{*}(X)  \stackrel{\cong}{\longrightarrow} 
  bu_{*}(\Sigma^{\infty} BS^{1}  \wedge X), \\
  \\
       bu_{*}(\Sigma^{\infty} W_{p}) \otimes_{{\mathbb Z}[u] }  bu_{*}(X)  \stackrel{\cong}{\longrightarrow}   bu_{*}(\Sigma^{\infty} W_{p}  \wedge X)  .
 \end{array}  \]

Smashing the cofibration of \S\ref{2.1} with $X$ and applying the argument of \cite{At62} yields the following K\"{u}nneth formula:
\end{em}
\end{dummy}
\begin{theorem}{$_{}$}
\label{2.8}
\begin{em}

There is a natural short exact sequence
\[   \begin{array}{l} 
0 \longrightarrow  bu_{*}(\Sigma^{\infty} B{\mathbb Z}/p) \otimes_{{\mathbb Z}[u] }  bu_{*}(X)  \longrightarrow     bu_{*}(\Sigma^{\infty} B{\mathbb Z}/p) \wedge X)    \\
\\
\hspace{100pt}  \longrightarrow  
 {\rm Tor}_{{\mathbb Z}[u]}^{1}( bu_{*}(\Sigma^{\infty} B{\mathbb Z}/p) , bu_{*}(X) )[1]    \longrightarrow  0   
  \end{array}    \]
as well as similar exact sequences for $bu{\mathbb Z}_{p}$ and $lu$.
\end{em}
\end{theorem}
\begin{example}
\label{2.9}
\begin{em}

Let $p$ be a prime. As in Corollary \ref{2.5}, let $lu_{*}$ be as in \S\ref{2.4} when $p$ is odd or $lu = (bu{\mathbb Z}_{2})_{*}$ when $p=2$. In Theorem \ref{2.8} set $X = \Sigma^{\infty} B{\mathbb Z}/p$. Then $  lu_{2*}(\Sigma^{\infty} (B{\mathbb Z}/p \wedge B{\mathbb Z}/p)) $ comes entirely from the left-hand graded group which is generated by $\{ v_{2i-1}  \otimes v_{2j-1} \ i, j \geq 1\}$ but if $2i-1  > 2(p-1)$ then 
\[   \begin{array}{ll}
pv_{2i-1}  \otimes v_{2j-1} & = v v_{2i-1 - 2(p-1)}  \otimes v_{2j-1} \\
& = v_{2i-1 - 2(p-1)}  \otimes v v_{2j-1}  \\
&=  p v_{2i-1 - 2(p-1)}   \otimes v_{2(p-1) + 2j-1} 
\end{array} \] 
which is zero by induction and similarly if $2j-1  > 2(p-1)$. Therefore
$  lu_{2m}(\Sigma^{\infty} (B{\mathbb Z}/p \wedge B{\mathbb Z}/p)) $ is the graded ${\mathbb F}_{p}$-vector space spanned by $v_{1} \otimes v_{2m-1}, \ldots , v_{2m-1} \otimes v_{1}$ which are linearly independent, being detected by the canonical homomorphism to $ H_{2m}(\Sigma^{\infty} (B{\mathbb Z}/p \wedge B{\mathbb Z}/p); {\mathbb Z}/p) $. Therefore for each $m \geq 1$ 
 \[  {\rm dim}_{{\mathbb F}_{p}}(lu_{2m}(\Sigma^{\infty} (B{\mathbb Z}/p \wedge B{\mathbb Z}/p)) ) = 
 {\rm dim}_{{\mathbb F}_{p}}( \pi_{2m}(    \vee_{i,j > 0} \   \Sigma^{2i+2j-2} H{\mathbb Z}/p ) ) . \]
 
 Similarly  $  lu_{2*+1}(\Sigma^{\infty} (B{\mathbb Z}/p \wedge B{\mathbb Z}/p)) $ comes entirely from the right-hand graded group 
 \[ T_{*} =  {\rm Tor}_{{\mathbb Z}_{p}[v]}^{1}( lu_{*}(\Sigma^{\infty} B{\mathbb Z}/p) , 
  lu_{*}(\Sigma^{\infty} B{\mathbb Z}/p) )[1]   . \]  
  
For $1 \leq 2i-1 \leq 2p-3$ let $Y_{i}$ denote the ${\mathbb Z}_{p}[v]$-submodule of $lu_{*}(\Sigma^{\infty} B{\mathbb Z}/p)$ generated by $\{ v_{2k(p-1)+2i-1} \ | \ k \geq 0\}$ so that $lu_{*}(\Sigma^{\infty} B{\mathbb Z}/p) \cong \oplus_{i} Y_{i}$
with a corresponding decomposition $T_{*} \cong  \oplus_{i} T_{i , *}$ .
This decomposition has a well-known geometric origin (\cite{Os89} \S2).

A free ${\mathbb Z}_{p}[v]$-module resolution is given by
\[   0 \longrightarrow  \oplus_{j=0}^{\infty} \  {\mathbb Z}_{p}[v] \langle a_{j}   \rangle   \stackrel{d}{\longrightarrow }  
 \oplus_{j=0}^{\infty} \  {\mathbb Z}_{p}[v] \langle b_{j}   \rangle  \stackrel{\epsilon}{\longrightarrow } 
Y_{i}  \longrightarrow   0   \]
 where $a_{j}  , b_{j} $ have internal degree $2j(p-1) + 2i-1$, $\epsilon(b_{j}) = v_{2j(p-1)+2i-1}, d(a_{0}) = pb_{0}$ and $d(a_{j}) = pb_{j} - v b_{j-1}$ for $j \geq 1$. Therefore
\[ T_{i, *} =  {\rm Ker}( 1 \otimes d :    \oplus_{j=0}^{\infty} \   lu_{*}(\Sigma^{\infty} B{\mathbb Z}/p)\langle a_{j}   \rangle  \longrightarrow    \oplus_{j=0}^{\infty} \   lu_{*}(\Sigma^{\infty} B{\mathbb Z}/p)\langle b_{j}   \rangle    ) . \] 
For $2m \geq 2i-1$ write $2m-2i+1 = 2t(p-1) + 2j-1$ with $1 \leq 2j-1 \leq 2p-3$. Then
\[   T_{i, 2m} =  {\mathbb Z}/p^{t+1} \langle   v_{2t(p-1) + 2j-1} a_{0} + v_{2(t-1)(p-1) + 2j-1} a_{1} + \ldots + v_{2j-1} a_{t} \rangle .  \]
From this, for $1 \leq i \leq p-1$, one finds that
\[    {\rm Tor}_{{\mathbb Z}_{p}[v]}^{1}( lu_{*}(\Sigma^{\infty} B{\mathbb Z}/p) , 
 Y_{i} )[1]  \cong  lu_{2 * +1}(\Sigma^{2i} B{\mathbb Z}/p)   \]
 and therefore
\[    lu_{2 * +1}(\Sigma^{\infty} B{\mathbb Z}/p \wedge B{\mathbb Z}/p)   \cong  \oplus_{i=1}^{p-1}  \  lu_{2 * +1}(\Sigma^{2i} B{\mathbb Z}/p) .\]
Adding together the suspensions of $lu$ as in \S\ref{2.4} yields a similar isomorphism for $bu{\mathbb Z}_{p}$ and therefore there is a ${\mathbb Z}[u]$-module isomorphism
\[   bu_{2*+1}(\Sigma^{\infty} (B{\mathbb Z}/p \wedge B{\mathbb Z}/p)) \cong  
 \oplus_{i=1}^{p-1}  \  bu_{2*+1}(\Sigma^{\infty} (S^{2i} \wedge B{\mathbb Z}/p))     \]
 and $ \oplus_{i=1}^{p-1}  bu_{2*}(\Sigma^{\infty} (S^{2i} \wedge B{\mathbb Z}/p))  = 0$.
 
 From Example \ref{2.9} we have an isomorphism
  \[ lu_{2*}(\Sigma^{\infty} (B{\mathbb Z}/p \wedge B{\mathbb Z}/p))  \cong 
 \pi_{2*}(    \vee_{i,j > 0} \   \Sigma^{2i+2j-2} H{\mathbb Z}/p )   \] 
 and therefore
  \[      bu_{2*}( \Sigma^{\infty} (B{\mathbb Z}/p \wedge B{\mathbb Z}/p)
\cong    \pi_{2*}(  
  \vee_{a= 0}^{p-2} \vee_{i,j > 0} \   \Sigma^{2a+2i+2j-2} H{\mathbb Z}/p )    . \]
\end{em}
\end{example}
\begin{lemma}{$_{}$}
\label{2.10}
\begin{em}

The homomorphism induced by the multiplication $\mu$ in the group ${\mathbb Z}/p$ is injective
\[   \mu_{*} :  bu_{2m+1}(\Sigma^{\infty} (B{\mathbb Z}/p \wedge B{\mathbb Z}/p)) \longrightarrow  bu_{2m+1}(\Sigma^{\infty} (B{\mathbb Z}/p)) . \]
\end{em}
\end{lemma}
\vspace{2pt}

{\bf Proof}
\vspace{2pt}

For simplicity we prove this only for $p=2$. The proof, which uses $KU$, may be modified for odd primes but requires a more careful analysis of the splittings of \S\ref{2.4} and (\cite{Os89} \S2) in relation to the embedding
\[ KU_{2m+1}(\Sigma^{\infty} (B{\mathbb Z}/p \wedge B{\mathbb Z}/p)) \subset  {\rm Hom}( R({\mathbb Z}/p \times {\mathbb Z}/p) , {\mathbb Z}/p^{\infty})  \]

 By Example \ref{2.8}, multiplication by $u$ is injective in odd dimensions so that
\[  bu_{2m+1}(\Sigma^{\infty} (B{\mathbb Z}/p \wedge B{\mathbb Z}/p)) \longrightarrow  KU_{2m+1}(\Sigma^{\infty} (B{\mathbb Z}/p \wedge B{\mathbb Z}/p)) \]
is injective, because it is localisation by inverting $u$. We shall use this observation to show that
\[   \mu_{*} : KU_{2m+1}(\Sigma^{\infty} (B{\mathbb Z}/2 \wedge B{\mathbb Z}/2)) \longrightarrow  KU_{2m+1}(\Sigma^{\infty} (B{\mathbb Z}/2 )) \]
is injective, which suffices to prove the result when $p=2$. Since  $B{\mathbb Z}/2 = {\mathbb RP}^{\infty}$ a skeletal approximation to the multiplication gives
\[   \mu:   \Sigma^{\infty}({\mathbb RP}^{2r}   \wedge  {\mathbb RP}^{2v})  \longrightarrow \Sigma^{\infty}{\mathbb RP}^{2r+2v} .  \]
Consider the effect on reduced, periodic complex K-theory
\[   \mu^{*} : \tilde{KU}^{0}({\mathbb RP}^{2r+2v})   \cong  {\mathbb Z}/2^{r+v}  \longrightarrow
\tilde{KU}^{0}(  {\mathbb RP}^{2r}   \wedge  {\mathbb RP}^{2v} ) \cong  {\mathbb Z}/2^{{\rm min}(r, v)}   . \]
If $L$ is the Hopf line bundle then $\mu^{*}(L-1) = (L-1) \otimes (L-1)$ so that $\mu^{*}$ is onto and, by the universal coefficient formula for $KU_{*}, KU^{*}$,
\[   \mu_{*} :  \tilde{KU}_{2m+1}(  {\mathbb RP}^{2r}   \wedge  {\mathbb RP}^{2v} ) \cong  {\mathbb Z}/2^{{\rm min}(r, v)}
\longrightarrow  \tilde{KU}_{2m+1}({\mathbb RP}^{2r+2v})   \cong  {\mathbb Z}/2^{r+v}  \]
is injective. Letting $r, s$ tend to infinity yields the result.  $ \Box $
\begin{dummy}
\label{2.11}
\begin{em}

We have a cofibration of spectra $\Sigma^{2}bu \longrightarrow bu \longrightarrow  H{\mathbb Z}$. By Example \ref{2.9} and Lemma \ref{2.10} the composition
\[     bu \wedge \Sigma^{\infty} (B{\mathbb Z}/p \wedge B{\mathbb Z}/p) \stackrel{1 \wedge \mu}{\longrightarrow} bu \wedge \Sigma^{\infty} (B{\mathbb Z}/p )  \stackrel{k}{ \longrightarrow}   H{\mathbb Z} \wedge \Sigma^{\infty} (B{\mathbb Z}/p)   \]
is trivial on homotopy groups. Therefore, when $p=2$,  $(1 \wedge \mu)_{*}$ induces an isomorphism 
\[   \tilde{\mu}_{*} :  bu_{2*+1}(\Sigma^{\infty} (B{\mathbb Z}/2 \wedge B{\mathbb Z}/2)) 
\stackrel{\cong}{\longrightarrow}  bu_{2*+1}(\Sigma^{\infty} (S^{2} \wedge  B{\mathbb Z}/2)) .  \]
Similarly at odd primes, using the multiplication $\mu$ together with the stable homotopy splittings of $\Sigma^{\infty} B{\mathbb Z}/p$ \cite{Os89}, yields an isomorphism
\[   \tilde{\mu}_{*} :  bu_{2*+1}(\Sigma^{\infty} (B{\mathbb Z}/p \wedge B{\mathbb Z}/p)) 
\stackrel{\cong}{\longrightarrow}  \oplus_{i=1}^{p-1} \  bu_{2*+1}(\Sigma^{\infty} (S^{2i} \wedge  B{\mathbb Z}/p)) .  \]

By Example \ref{2.9}, the ${\mathbb F}_{p}$-vector space $bu_{2*}(\Sigma^{\infty} (B{\mathbb Z}/p \wedge B{\mathbb Z}/p))$ is detected in mod $p$ homology and there is a map of spectra
\[  h:      bu \wedge \Sigma^{\infty} (B{\mathbb Z}/p \wedge B{\mathbb Z}/p) \longrightarrow  
  (\vee_{a= 0}^{p-2} \  \vee_{i,j > 0} \   \Sigma^{2a+2i+2j-2} H{\mathbb Z}/p )   \]
 which induces an isomorphism on even dimensional homotopy. Therefore we obtain the following result:
\end{em}
\end{dummy}
\begin{theorem}{(\cite{Os89}; see also \cite{Br99} and \cite{JW97})}
\label{2.12}
\begin{em}

 There is an isomorphism
 \[  \begin{array}{l}
 ( \tilde{\mu}_{*}, h_{*}) :     bu_{*}( \Sigma^{\infty} (B{\mathbb Z}/p \wedge B{\mathbb Z}/p))
\stackrel{\cong}{  \longrightarrow }   \\
\\
\hspace{60pt}    \oplus_{i=1}^{p-1} \  bu_{*}(   \Sigma^{\infty} (S^{2i} \wedge  B{\mathbb Z}/p))  \oplus 
  \oplus_{a= 0}^{p-2} \oplus_{i,j > 0} \   \pi_{*}(\Sigma^{2a+2i+2j-2} H{\mathbb Z}/p  )   . 
  \end{array}  \]
\end{em}
\end{theorem}
\begin{remark}
\label{2.13}
\begin{em}

The composition of maps of spectra $k(1 \wedge \mu)$
used in \S\ref{2.11} is not nullhomotopic, although is it zero on homotopy groups. It is for this reason that our method does not yield a homomorphism $  \tilde{\mu}_{*} $ induced by a map of spectra.
\end{em}
\end{remark}

\section{The connective orthogonal case}
\begin{dummy}
\label{3.1}
\begin{em}

In this section we shall concentrate on $p=2$ and connective orthogonal K-theory $bo$. Consider the following commutative diagram of spectra of horizontal and vertical cofibrations in which $c$ is complexification and $\eta$ is multiplication by the generator of $\pi_{1}(bo)$. The notation for $bo \langle 1 \rangle$ is taken from \cite{Br99}.
\newline
{\tt    \setlength{\unitlength}{0.92pt}
\begin{picture}(459,186)
\thinlines   
              \put(305,45){\vector(1,0){70}}
              \put(88,45){\vector(1,0){67}}
              \put(398,45){$H{\mathbb Z}$}
              \put(230,45){$bo$}
              \put(9,45){$bo \langle 1 \rangle$}
               \put(9,-45){$\Sigma^{2}bu$}
              \put(230,-45){$bu$}
              \put(238,124){\vector(0,-1){45}}
              \put(20,121){\vector(0,-1){45}}
              \put(88,165){\vector(1,0){63}}
              \put(232,161){$\Sigma bo$}
              \put(15,161){$\Sigma bo$}
              \put(404,27){\vector(0,-1){46}}
              \put(238,24){\vector(0,-1){45}}
              \put(20,21){\vector(0,-1){45}}
                \put(305,-43){\vector(1,0){70}}
              \put(88,-43){\vector(1,0){67}}
               \put(398,-45){$H{\mathbb Z}$}
               \put(230,0){$c$}
               \put(396,0){$1$}
               \put(25,0){$\tilde{c}$}
               \put(230,100){$\eta$}
               \put(25,100){$\tilde{\eta}$}
               \put(110,155){$1$}
\end{picture}}
\newline
\vspace{40pt}

We have the following table of (reduced) orthogonal connective K-theory groups:
\[ \begin{array}{||c|c||}
\hline
bo_{8n}({\mathbb RP}^{\infty} ) & 0  \\
\hline
bo_{8n+1}({\mathbb RP}^{\infty} ) & {\mathbb Z}/2   \\
\hline
bo_{8n+2}({\mathbb RP}^{\infty} ) & {\mathbb Z}/2   \\
\hline
bo_{8n+3}({\mathbb RP}^{\infty} ) & {\mathbb Z}/2^{4n+3}   \\
\hline
bo_{8n+4}({\mathbb RP}^{\infty} ) & 0  \\
\hline
bo_{8n+5}({\mathbb RP}^{\infty} ) & 0  \\
\hline
bo_{8n+6}({\mathbb RP}^{\infty} ) & 0  \\
\hline
bo_{8n+7}({\mathbb RP}^{\infty} ) & {\mathbb Z}/2^{4n+4}   \\
\hline
\end{array}  \]
The graded group $bo_{*}({\mathbb RP}^{\infty} ) $ is a module over 
\[  bo_{*}(S^{0}) = {\mathbb Z}[\eta,\alpha,\beta]/(2 \eta, \eta^3, \eta \alpha, \alpha^2 - 4 \beta)  \]
 and multiplication by $\eta$ is nontrivial from dimension $8n+1$ to $8n+2$ and from $8n+2$ to $8n+3$. Multiplication by $\alpha$ has kernel of order $4$ from dimension $8n+3$ to $8n+7$ and is 
one-one from dimension $8n+7$ to $8n+11$. Multiplication by 
$\beta$ is always one-one.

The central horizontal cofibration yields a long exact sequence of reduced homology theories
\[   \ldots \longrightarrow   bo \langle 1 \rangle_{i}({\mathbb RP}^{\infty})  \longrightarrow   bo_{i}({\mathbb RP}^{\infty}) \longrightarrow   H_{i}( {\mathbb RP}^{\infty} ; {\mathbb Z} )   \longrightarrow   \ldots \]
and there is a factorisation  $bo \longrightarrow  bu \longrightarrow  H{\mathbb Z}$. Using the fact that $H_{i}( {\mathbb RP}^{\infty} ; {\mathbb Z} )  \cong {\mathbb Z}/2$ for odd $i > 0$ and is zero otherwise we may calculate $ bo \langle 1 \rangle_{*}({\mathbb RP}^{\infty}) $. In addition we may double-check the results from the long exact homotopy sequence of the left-hand vertical fibration in the diagram of \S\ref{3.1}
\[   \ldots  \longrightarrow  bo_{i-1}({\mathbb RP}^{\infty}) \longrightarrow\longrightarrow   bo \langle 1 \rangle_{i}({\mathbb RP}^{\infty})  \longrightarrow   bu_{i-2}({\mathbb RP}^{\infty}) \longrightarrow     \ldots   . \]
Diagram chasing yields the following table:
\[  \begin{array}{||c|c||}
\hline
bo<1>_{8n}({\mathbb RP}^{\infty} )  \  n \geq 0 &  {\mathbb Z}/2 \cong H_{8n+1}({\mathbb RP}^{\infty} ; {\mathbb Z}) \stackrel{\cong}{\longrightarrow}  bo<1>_{8n}({\mathbb RP}^{\infty} )\\
\hline
bo<1>_{1}({\mathbb RP}^{\infty} )   &  0  \\
\hline
bo<1>_{8n+1}({\mathbb RP}^{\infty} )   \  n \geq 1   & {\mathbb Z}/2  \cong bo<1>_{8n+1}({\mathbb RP}^{\infty} )  
\stackrel{\cong}{\longrightarrow}  bo_{8n+1}({\mathbb RP}^{\infty} ) \\
\hline
bo<1>_{8n+2}({\mathbb RP}^{\infty} ) &  {\mathbb Z}/2  \cong bo<1>_{8n+2}({\mathbb RP}^{\infty} )  
\stackrel{\cong}{\longrightarrow}  bo_{8n+2}({\mathbb RP}^{\infty} )   \\
\hline
bo<1>_{8n+3}({\mathbb RP}^{\infty} ) & {\mathbb Z}/2^{4n+2} \cong  bo<1>_{8n+3}({\mathbb RP}^{\infty} )  
\stackrel{2(2s+1)}{\longrightarrow}   bu_{8n+3}({\mathbb RP}^{\infty} )  \\
\hline
bo<1>_{8n+4}({\mathbb RP}^{\infty} ) &  {\mathbb Z}/2 \cong H_{8n+5}({\mathbb RP}^{\infty} ; {\mathbb Z}) \stackrel{\cong}{\longrightarrow}  bo<1>_{8n+4}({\mathbb RP}^{\infty} ) \\
\hline
bo<1>_{8n+5}({\mathbb RP}^{\infty} ) & 0  \\
\hline
bo<1>_{8n+6}({\mathbb RP}^{\infty} ) & 0  \\
\hline
bo<1>_{8n+7}({\mathbb RP}^{\infty} ) & {\mathbb Z}/2^{4n+3}  \cong  bo<1>_{8n+7}({\mathbb RP}^{\infty} )  
\stackrel{1-1}{\longrightarrow}   bo_{8n+7}({\mathbb RP}^{\infty} )  \\
\hline
\end{array}  \]
We can now state the main result of this section (see also \cite{Br99}), whose proof will be sketched in \S\ref{3.7}.
\end{em}
\end{dummy}
\begin{theorem}{$_{}$}
\label{3.2}
\begin{em}

There are homomorphisms of graded groups
 \[  \tilde{\mu}_{*  } : bo_{*}({\mathbb RP}^{\infty}  \wedge {\mathbb RP}^{\infty} ) \longrightarrow  bo \langle 1 \rangle_{*}({\mathbb RP}^{\infty} ),  \]
characterised in \S\ref{3.7}, and 
\[  \tilde{h}_{*} : bo_{m}({\mathbb RP}^{\infty}  \wedge {\mathbb RP}^{\infty}) \longrightarrow  \pi_{*}( \vee_{i,j > 0} \   \Sigma^{2i+4j-2} H{\mathbb Z}/2 )  \]
such that
  \[  (\tilde{\mu}_{*} ,  \tilde{h}_{*} ) :   bo_{*}({\mathbb RP}^{\infty}  \wedge {\mathbb RP}^{\infty})  \longrightarrow  bo \langle 1 \rangle_{*}({\mathbb RP}^{\infty}  ) \oplus   \pi_{*}( \vee_{i,j > 0} \   \Sigma^{2i+4j-2} H{\mathbb Z}/2 )  \]
  is an isomorphism. 
\end{em}
\end{theorem}
\begin{dummy}{The Bott sequence versus mod $2$ homology}
\label{3.4}
\begin{em}

The subalgebra ${\mathcal B}$ generated in the mod $2$ Steenrod algebra by $Sq^{1}$ and $Sq^{2}$ has dimension eight and contained the  exterior subalgebra ${\mathcal E}=E(Sq^{1}, Sq^{0,1}) = \{ 1 , Sq^{1},   Sq^{1} Sq^{2}+ Sq^{2}Sq^{1} , Sq^{2}Sq^{2} \}$.

Consider the Bott sequence
\[  \cdots \longrightarrow  bo_{i}(X)  \stackrel{c}{\longrightarrow}  bu_{i}(X) \longrightarrow  bo_{i-2}(X)
\stackrel{\eta_{*}}{\longrightarrow}  bo_{i-1}(X)  \longrightarrow  \cdots  \]
which is isomorphic to the homotopy sequence of the cofibration
\[    bo \wedge X  \longrightarrow  bo \wedge \Sigma^{-2} {\mathbb CP}^{2} \wedge X  \longrightarrow  bo \wedge S^{2} \wedge X  , \]
where the middle spectrum is identified with $bu \wedge X$ via an equivalence due to Anderson-Wood \cite{Sn09}.

The following commutative diagram is easy to establish.
\newline
\begin{picture}(459,186)
              \put(118,45){\vector(1,0){40}}
              \put(200,-45){$ {\rm Hom}_{\mathcal{B}}( H^{i-2}( X; {\mathbb Z}/2)  ,   {\mathbb Z}/2  ) $}
              \put(9,-45){$bo_{i-2}( X )$}
               \put(9,45){$ bu_{i}(X)$}
              \put(200,45){$ {\rm Hom}_{\mathcal{E}}( H^{i}( X; {\mathbb Z}/2)  ,   {\mathbb Z}/2  ) $}
              \put(238,124){\vector(0,-1){45}}
              \put(20,121){\vector(0,-1){45}}
              \put(88,165){\vector(1,0){63}}
              \put(200,161){$  {\rm Hom}_{\mathcal{B}}( H^{i}(X; {\mathbb Z}/2)  ,   {\mathbb Z}/2  ) $}
              \put(15,161){$ bo_{i}(X)$}
              \put(238,24){\vector(0,-1){45}}
              \put(20,21){\vector(0,-1){45}}
              \put(88,-43){\vector(1,0){67}}
               \put(226,0){$\tilde{\lambda}$}
               \put(226,100){$\tilde{\phi}$}
               \put(25,100){$c$}
\end{picture}
\newline
\vspace{40pt}

The horizontal maps are induced by the canonical map $bo \longrightarrow  H{\mathbb Z}/2$, $\tilde{\phi}(h) = h$ and $\tilde{\lambda}(g)(x) = g( Sq^{2}(x))$.

By Theorem \ref{2.12} we know that the middle horizontal map is an isomorphism when $i$ is even. The following result is straightforward.
\end{em}
\end{dummy}
\begin{prop}{$_{}$}
\label{3.5}
\begin{em}

(i)  \  When $X = {\mathbb RP}^{\infty}  \wedge {\mathbb RP}^{\infty}$ the sequence
\[  \begin{array}{l}
 0   \longrightarrow  {\rm Hom}_{\mathcal{B}}( H^{i}(X ; {\mathbb Z}/2)  ,   {\mathbb Z}/2  )  
 \stackrel{\tilde{\phi}}{  \longrightarrow }  {\rm Hom}_{\mathcal{E}}( H^{i}(X; {\mathbb Z}/2)  ,   {\mathbb Z}/2  )  \\
 \\
 \hspace{65pt}   \stackrel{\tilde{\lambda}}{\longrightarrow } {\rm Hom}_{\mathcal{B}}( H^{i-2}(X; {\mathbb Z}/2)  ,   {\mathbb Z}/2  ) 
 \longrightarrow  0 .  
\end{array}  \]
is exact.

(ii)  \   For $i \geq 2$ 
\[   {\rm dim}_{{\mathbb F}_{2}}( {\rm Hom}_{\mathcal{B}}( H^{2i}( {\mathbb RP}^{\infty}  \wedge {\mathbb RP}^{\infty} ; {\mathbb Z}/2)  ,   {\mathbb Z}/2  )  ) = \left\{   
\begin{array}{l}
 (i-1)/2 \ {\rm if} \ i  \  {\rm is \ odd}, \\
\\
1 + (i/2)  \ {\rm if} \ i  \   {\rm is \ even}
\end{array} \right. \]

(iii)  \   For $i \geq 2$
\[   {\rm dim}_{{\mathbb F}_{2}}( {\rm Hom}_{\mathcal{E}}( H^{2i}( {\mathbb RP}^{\infty}  \wedge {\mathbb RP}^{\infty} ; {\mathbb Z}/2)  ,   {\mathbb Z}/2  )  ) = i . \]
\end{em}
\end{prop}

We shall also need the following result.
\begin{prop}{$_{}$}
\label{3.6}
\begin{em}

Define 
\[   X_{n} = \# \{ i,j \geq 1  \  |  \  2i-1+4j-1 = 2n  \}  .\]
Then
\[ X_{n} = \left\{  \begin{array}{ll}
n/2 & {\rm if} \ n \ {\rm is \ even},  \\
\\
(n-1)/2  & {\rm if} \ n \geq 1 \ {\rm is \ odd},
\end{array} \right.  \]
\end{em}
\end{prop}
\begin{dummy}{Sketch proof of Theorem \ref{3.2}}
\label{3.7}
\begin{em}

The mod $2$ Hurewicz homomorphism induces a homomorphism
\[ h_{*}  :  bo_{m}({\mathbb RP}^{\infty} \wedge {\mathbb RP}^{\infty})  \longrightarrow   {\rm Hom}_{\mathcal{B}}( H^{2i}( {\mathbb RP}^{\infty}  \wedge {\mathbb RP}^{\infty} ; {\mathbb Z}/2)  ,   {\mathbb Z}/2  ) .\]
The results in the central columns of the following table are proved by induction on dimension using Theorem \ref{2.12}, the map between Bott sequences induces by
\[ \mu : \Sigma^{\infty}  {\mathbb RP}^{\infty}  \wedge {\mathbb RP}^{\infty} \longrightarrow  \Sigma^{\infty}  {\mathbb RP}^{\infty}  \]
and the results of \S\ref{3.4}, Proposition \ref{3.5} and Proposition \ref{3.6} concerning mod $2$ homology.
\[ \begin{array}{|c|c|c|c|}
\hline
m & bo_{m}({\mathbb RP}^{\infty} \wedge {\mathbb RP}^{\infty}) &  h_{*}  & 
 bo \langle 1 \rangle_{m}({\mathbb RP}^{\infty}) \\
\hline
\hline
2 & {\mathbb Z}/2&{\rm isom}&  {\mathbb Z}/2\\
\hline
3 \leq 8n+3 & {\mathbb Z}/2^{4n+3}  & - & {\mathbb Z}/2^{4n+3}\\
\hline
4 \leq 8n+4 & ({\mathbb Z}/2)^{2n+2} &{\rm isom}& {\mathbb Z}/2 \\
\hline
 5 \leq 8n+ 5& 0 & -  & 0  \\
\hline 
6 \leq 8n+ 6 &  ({\mathbb Z}/2)^{2n+1} &{\rm isom}& 0 \\
\hline 
7 \leq 8n+ 7 & {\mathbb Z}/2^{4n+4} & - &  {\mathbb Z}/2^{4n+4} \\
\hline 
8 \leq 8n & ({\mathbb Z}/2)^{2n+1}  &{\rm isom}& {\mathbb Z}/2 \\
\hline 
9 \leq 8n+ 1& {\mathbb Z}/2 & -  & {\mathbb Z}/2 \\
\hline 
10 \leq 8n+ 2 &({\mathbb Z}/2)^{2n+1}  &{\rm onto, \ not \ isom}& {\mathbb Z}/2 \\
\hline
\hline
\end{array}  \]
The group $bo_{m}({\mathbb RP}^{\infty} \wedge {\mathbb RP}^{\infty})  = 0$ for $m \leq 1$.

When $m = 8n+1, 8n+2,  8n+3, 8n+7$ it is important that the homomorphism $\tilde{\mu}_{m}$ is chosen so that the composition 
\[ bo_{m}({\mathbb RP}^{\infty}  \wedge {\mathbb RP}^{\infty} ) \longrightarrow  bo \langle 1 \rangle_{m}({\mathbb RP}^{\infty} )   \subseteq   bo_{m}({\mathbb RP}^{\infty} ) \]
is equal to $\mu_{m}$.
\end{em}
\end{dummy}

\end{document}